\documentclass[11pt]{amsart}
\usepackage{verbatim, latexsym, amssymb, amsmath,color}
\usepackage{epsfig} 
\usepackage{graphicx}
\usepackage{hyperref}

\title[The first min-max width of surfaces with boundary]{The first width of non-negatively curved surfaces with convex boundary}

\author{Sidney Donato}
\author{Rafael Montezuma}

\thanks{SD was supported by CAPES-Brasil 001. RM was supported by CNPq and Instituto Serrapilheira,
grant "New perspectives of the min-max theory for the area functional". \vspace{.1cm}\\  Instituto de Matem\'atica, Universidade Federal de Alagoas, Macei\'o, AL, 57072-970, Brazil. \textit{E-mail address:} sidney.silva@im.ufal.br \vspace{.1cm}\\ Universidade Federal do Cear\'{a}, Departamento de Matem\'{a}tica,  Fortaleza, CE 60455-760, Brazil. \textit{E-mail address:} montezuma@mat.ufc.br
}

\begin{document}

\maketitle

\begin{abstract}
In this paper, free boundary geodesic networks whose length realize the first min-max width of the length functional are investigated. This functional acts on the space of relative flat $1$-dimensional cycles modulo $2$ in a compact surface with boundary. The widths are special critical values of the volume functional in some class of submanifolds which naturally arise in the Min-max Theory of Almgren and Pitts. 

The main result of this work concerns the existence of a geodesic network with a rather simple structure which realizes the first width of a surface with non-negative sectional curvature and strictly convex boundary. More precisely, it is either a simple geodesic meeting the boundary orthogonally, or a geodesic loop with vertex at a boundary point determining two equal angles with that boundary curve.  
 \end{abstract}
\vspace{.1cm}

\section{Introduction}

The min-max theory of Almgren and Pitts is a Morse theory for the $k$-dimensional volume defined in a space of generalized submanifolds of a Riemannian manifold. This theory was developed to address foundational questions about the existence of critical points of this functional: minimal submanifolds. Birkhoff was the first to consider min-max methods in the setting of the length functional $L$, $k=1$, to prove that every Riemannian two-sphere contains at least one closed geodesic. He considered the quantity 
$$
L = \inf_{\Phi}\left( \sup \{L(\Phi(t)): t \in [0,1]\}\right),
$$
where the infimum is taken over the set of all continuous families $\Phi = \{\Phi(t)\}$, $0\leq t \leq 1$, of closed curves $\Phi(t) \subset S^2$ contained in the two-sphere, such that $\Phi(0)$ and $\Phi(1)$ are constant curves and $\Phi$ is not homotopically trivial. It happens that the quantity $L$ is realized by the length of a closed geodesic.

Almgren \cite{almgren-varifolds} developed the theory showing that every closed Riemannian manifold contains at least one stationary integral $k$-dimensional varifold, for each $k$ between $1$ and $n$, where $dim(M) = n+1$. Pitts \cite{pitts} and Schoen and Simon \cite{SchSim} improved the regularity in the case of codimension one and $n+1\geq 3$, allowing singular sets of codimension $7$. These hypersurfaces realize a min-max quantity. More precisely, there exists a sequence of min-max numbers $\{\omega_p(M)\}$, for $p \in \mathbb{N}^{\ast}$, which is a nonlinear analog of the spectrum of its Laplace-Beltrami operator. These numbers are the so called $p$-widths, and they are realized by the volume of minimal hypersurfaces.

The relevance of the $p$-widths is shown in connection with results on the existence of infinitely many critical points of the volume functional, and related to spectral properties. See \cite{ChoMan-surf}, \cite{Don}, \cite{KeiMarNev}, \cite{LioMarNev}, \cite{marques-neves-infinitely}, \cite{MarNevSon} and references therein for further properties of the min-max $p$-widths.

In the case of dimension $n+1 = 2$, Pitts \cite{pitts-dim1} proved that the relevant min-max invariant is realized by the length of a geodesic network. Geodesic networks generalize closed geodesics, allowing for finite unions of geodesic segments, see Section \ref{sect-struct}. Other interesting proofs of Pitts' result and developments appeared in the works of Aiex \cite{aiex}, Calabi and Cao \cite{CalCao}, Zhou and Zhu \cite{ZhoZhu}, and references therein. The article \cite{ZhoZhu} and the article of Ketover and Liokumovich \cite{KetLio} include also the case of curves with constant geodesic curvature different from zero, and Cheng and Zhou \cite{DaRZho} and Asselle and Benedetti \cite{AssBen} considered the case of prescribed geodesic curvature. 

Different approaches to the Almgren-Pitts theory have appeared. One particularly successful approach uses a phase-transition regularization of the volume functional, see \cite{ChoMan}, \cite{GasGua}, \cite{GasGua2}, \cite{Gua}, and references therein. A sequence of $p$-widths is also introduced in this setting. Mantoulidis \cite{Man}, for $p=1$, and Chodosh and Mantoulidis \cite{ChoMan-surf} proved that these $p$-widths are realized by unions of closed immersed geodesics, possibly with multiplicities. Dey \cite{Dey} showed that the $p$-widths of the phase-transition approach agree with the ones considered in \cite{LioMarNev}, \cite{marques-neves-infinitely}. Therefore, \cite{Man} and \cite{ChoMan-surf} improve and extend some of the results obtained in the works \cite{pitts-dim1}, \cite{aiex}, \cite{CalCao}, \cite{ZhoZhu}.

In the case of compact manifolds $M$ with boundary, Li and Zhou \cite{LiZhou} extended the theory by showing that a min-max quantity defined in the space of relative cycles is realized by the volume of an embedded free boundary minimal hypersurface; a minimal hypersurface $\Sigma$ whose boundary $\partial \Sigma$ is contained in $\partial M$, and such that $\Sigma$ intersects $\partial M$ orthogonally along $\partial \Sigma$.

The first author investigated \cite{Don} the free boundary setting of the Almgren-Pitts theory for curves. The space of curves that is considered in this setting is $\mathcal{Z}_{1,rel}(M, \partial M; \mathbb{Z}_2)$, the space of relative flat cycles modulo $2$. Almgren \cite{almgren-groups} showed that this space is weakly homotopy equivalent to $\mathbb{R}P^{\infty}$. Let $\overline{\lambda}$ be the generator of the cohomology ring of that space of curves. The $p$-sweepouts are continuous maps $\Phi : X \rightarrow \mathcal{Z}_{1,rel}(M, \partial M; \mathbb{Z}_2)$ such that $p$-th cup power of $\overline{\lambda}$ satisfies that $\Phi^{\ast}(\overline{\lambda})^p$ is non-zero in $H^p(X; \mathbb{Z}_2)$. The $p$-width of $M$, denoted by $\omega_p(M)$, is defined as the infimum of numbers $\ell$ such that there exists a $p$-sweepout by curves of length at most $\ell$. Let us write $\omega(M) = \omega_1(M)$.

In \cite{Don}, it is proven that, on compact surfaces $M$ with convex boundary, a free boundary stationary integral $1$-varifold $V$ which is $\mathbb{Z}_2$-almost minimizing in small annuli is a free boundary geodesic network. Hence, the $p$-widths can be approximated by lengths of free boundary geodesic networks.

Theorem \ref{thm-single-jct} of the present work implies that $\omega(M)$ is realized by the length of a free boundary geodesic network with at most one singular point. The main result of this paper is the following theorem.

\subsection{Theorem}\label{thm-main}
{\em Let $(M^2,g)$ be a compact Riemannian surface with non-empty boundary. Assume that the sectional curvature of $g$ is non-negative, and that its boundary is strictly convex. Then, there exists a free boundary geodesic network $V$ with at most one junction whose length is equal to $\omega(M)$. Moreover, $V$ has multiplicity one and the support of $V$ is either a geodesic loop with vertex at a point in $\partial M$, or a free boundary geodesic.}

A consequence of the Gauss-Bonnet Theorem is that $M$ is necessarily diffeomorphic to a closed disk. Next, examples showing that the two possible structures for the geodesic network of Theorem \ref{thm-main} hold are presented.

\begin{figure}[htp]
    \centering
    \includegraphics[width=6cm]{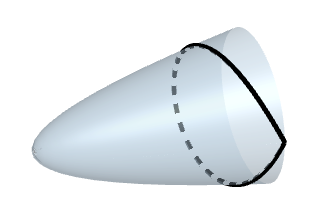}
    \caption{A positively curved convex surface for which the min-max width is realized by the length of a geodesic loop.}
    \label{Fig-1}
\end{figure}

Let $r$ be a smooth function defined on $[0, +\infty)$ such that: $r(0)=0$, $r\geq 0$, with right derivative $r^{\prime}_+(0) = +\infty$, $r^{\prime}(u)>0$, and $r^{\prime\prime}(u) \leq 0$, for all $u>0$. Assume that $r(u) = au+b$, for every $u\geq u_0$, where $a,b>0$. Let $M$ be obtained by the rotation of the graph of $r= r(u)$, $u \in [0, u_1]$, around the $u$-axis, for $u_1 > u_0$. The geodesic curvature of the boundary is positive, and the sectional curvature is non-negative. Consider $\alpha_{\theta} (u) = (u, r(u)\cos \theta, r(u) \sin \theta)$, and observe that the concatenations of $\alpha_{\theta}$ and $\alpha_{\theta+\pi}$ are precisely the free boundary geodesics of $M$. The family $\Phi_u(\theta) = \alpha_{\theta}(u)$ is a sweepout by curves of lengths bounded by $2\pi r(u_1)$. For $a$ small and $u_1$ large, $2\pi r(u_1)$ is lower than the length of those free boundary geodesics. Then, the free boundary geodesics do not realize the first width, and Theorem \ref{thm-main} implies that $\omega(M)$ coincides with the length of a geodesic loop, as in Figure \ref{Fig-1}.

In the cases that $M$ is a convex domain of the Euclidean plane or of the Euclidean two-sphere, the relative min-max width coincides with the length of a free boundary geodesic. These surfaces do not admit geodesic loops.

The method applied in the proof of Theorem \ref{thm-main} involves constructions of sweepouts by concatenations of sweepouts of certain convex regions. The curve shortening process of Birkhoff adapted to the boundary case plays a key role in these constructions. This process can also be applied to prove that, in the setting described in Theorem \ref{thm-main}, one has $\omega(M)< L(\partial M)$.

A different min-max approach was considered in Zhou \cite{Zhou-FB-Geod}, to show the existence of free boundary geodesics. He considered the min-max number associated to the length in the space of paths in $M$ with extremities in a fixed closed submanifold of general dimension and codimension. In particular, this number can be strictly larger than the $\omega(M)$ of Theorem \ref{thm-main}. See also related results of Gluck and Ziller \cite{GluZil}, Weinstein \cite{Wei}, and Nabutovsky and Rotman \cite{NabRot}. More recently, Li \cite{Li} devised a geometric flow which is the negative gradient flow for the length functional on the space of chords with end points lying on a fixed submanifold in Euclidean space. He gave simplified proofs of classical theorems on the existence of multiple free boundary geodesics. 

Besides the applications of min-max methods to prove existence results, they usually involve these interesting invariants which relate to other geometric quantities. In the setting of minimal surfaces in $3$-manifolds, the works of Marques and Neves \cite{MarNev-rigidity} and Ambrozio and the second author \cite{AmbMon} provide some examples of relation between the width and scalar curvature, or the width and the volume of the ambient. Liokumovich, Marques, and Neves \cite{LioMarNev} provided a general Weyl's law for the asymptotic behavior of the sequence of min-max widths. In general, there exists a deep relationship between properties of this sequence and the linear spectrum of the space. 

In the case of Riemannian metrics in the two-sphere, there are interesting results which relate its area and systole, the length of the shortest closed geodesic. See \cite{BeaRot}, \cite{CalCao}, and \cite{Cro} for instance. In positively curved two-spheres, the width coincides with the systole. Therefore, some of the results in \cite{BeaRot}, \cite{CalCao}, and \cite{Cro} can be interpreted as relationships between the first width and the area. Abbondandolo, Bramham, Hryniewicz and Salom\~ao \cite{ABHS} proved a sharp inequality involving the systole and the area of a pinched positively curved two-sphere. Their proofs involve techniques from both symplectic and Riemannian geometry. Adelstein and Vargas Pallete \cite{AdePal} controlled the systole on a two-sphere with non-negative curvature by three times the diameter. They also prove an optimal inequality involving systole and diameter in the setting of pinched metrics. Other relations between systole and diameter, including curvature free estimates, were obtained in Croke \cite{Cro}, Maeda \cite{Mae}, Nabutovsky and Rotman \cite{NabRot1}, Rotman \cite{Rot}, and Sabourau \cite{Sab}. 

In some arguments, such as in  \cite{ABHS} and \cite{CalCao}, it is relevant to observe that closed geodesic of shortest length must be simple. The present work gives a similar information about the curves which realize the relative Almgren-Pitts min-max width in non-negatively curved convex surfaces with boundary. 

Theorem \ref{thm-main} implies that the first min-max width of a triangle coincides with the length of its shortest height. Indeed, approximate the triangle $T$ by strictly convex planar regions $\Omega_i$ whose width is realized by the length of a free boundary segment $V_i$. These $V_i$ do not converge to a segment connecting points in the interior of two different sides of $T$. Moreover, they do not converge to a side of $T$ because, otherwise, the sweepouts by segments parallel to a fixed height of $T$ would have a lower maximum of lengths. The only remaining possibility is that $V_i$ converges to a height of $T$. In particular, it follows that the unit area triangle of highest $\omega(M)$ is the equilateral one.

The first author computed \cite{Don} the first four relative min-max widths of the flat disk of unit radius: $\omega_1(\mathbb{D}) = \omega_2(\mathbb{D}) = 2$ and $\omega_3(\mathbb{D}) = \omega_4(\mathbb{D}) = 4$. The Lusternik-Schnirelmann inequality and Theorem \ref{thm-main} together imply that $\omega_5(\mathbb{D})$ is strictly bigger than $\omega_4(\mathbb{D})$.

\textbf{Organization of the paper.} In Section \ref{sect-struct}, the structure of $1$-varifolds arising in the Min-max Theory of Almgren and Pitts is briefly revisited. Theorem \ref{thm-single-jct}  shows that there exists a free boundary stationary geodesic network with at most one junction realizing the width. In Section \ref{Sec-Birkhoff}, the version of the curve shortening process of Birkhoff adapted to the boundary case is applied in the construction of sweepouts of certain convex regions. In Section \ref{Sec-W-type}, some properties of general free boundary geodesic networks are established and applied in the construction of induced sweepouts of the entire surface. In Section \ref{sec-single-junction}, the structures of free boundary geodesic networks with at most one junction, and different from a geodesic loop and a free boundary geodesic are ruled out as being the geodesic network of Theorem \ref{thm-single-jct}. In Section \ref{sec-proof}, the proof of Theorem \ref{thm-main}, and of the inequalities $\omega(M)< L(\partial M)$ and $\omega_4(\mathbb{D}) < \omega_5(\mathbb{D})$ are presented.

\section{Free boundary geodesic networks arising in min-max theory}\label{sect-struct}

One of the main objects considered in this work is that of a geodesic network. The adapted version for the boundary case will play a key role. Let us start by recalling the definitions of these objects.

A $1$-varifold $V$ is said to be a geodesic network in an open subset $U$ of a closed Riemannian manifold $M$ if the restriction of $V$ to $U$ is stationary in the set $U$ and it can be written as the sum of finitely many varifolds induced by geodesic segments $\alpha_i$ with respective positive integer multiplicities $m_i$. The stationarity of $V$ means that it is a critical point of the length functional.

A geodesic network $V$ has a set of distinguished points, called junctions, which can be thought as the singular points of $V$. More precisely, a \textit{junction} is a point $p\in U$ such that $p \in \partial \alpha$ for some geodesic segment $\alpha$ of $V$. If the geodesic segments of $V$ arriving at $p$ can be paired in such a way that every pair has opposite directions and equal multiplicities, then the union of the two geodesics in each pair can be considered as a single geodesic in $V$, and $p$ is said to be a \textit{regular junction}. Besides of junctions, the geodesic segments may intersect at interior points, which are referred to as \textit{crossing points}.

The stationarity of $V$ in $U$ implies that at least three geodesics arrive at each junction. Moreover, at each junction $p$, the resultant vector $\sum_k m_{i_k} \vec{\alpha}_{i_k}$ of the unit vectors $\vec{\alpha}_{i_k}$ pointing in the direction of the geodesics of $V$ emanating from $p$, with respective multiplicities, is the zero vector.

The relevance of the concept of a geodesic network is emphasized in the work \cite{AllAlm} of Allard and Almgren, in which the structure of stationary integral $1$-varifolds is classified as that of the geodesic networks, see also Theorem 3.5 of \cite{aiex}. The following definition will be relevant in this work.

\subsection{Definition}\label{FBGN}
Let $M$ be a compact Riemannian manifold with non-empty boundary $\partial M$, and $U\subset M$ be a relatively open subset. A $1$-varifold $V$ is said to be a {\em free boundary geodesic network in} $U$ if $V$ is a sum of geodesic segments as above, and is stationary with respect to the variations supported in $U$ which preserve $(\partial M)\cap U$. At junctions $p \in U\cap (\partial M)$, the sum $\vec{v}_p = \sum_k m_{i_k} \vec{\alpha}_{i_k}$ of the unit vectors $\vec{\alpha}_{i_k}$ pointing inwards in the direction of the geodesics of $V$ emanating from $p$ is perpendicular to $\partial M$.    

The one-dimensional density $\theta^1(V,p)$ coincides with half of the number of geodesic segments of $V$ starting at $p$ (with multiplicities). The following regularity theorem combines results obtained in \cite{aiex}, \cite{AllAlm}, and \cite{Don}.

\subsection{Theorem}[Thms. 16 of \cite{Don}, 3.4 and 4.13 of \cite{aiex}]\label{structure-fb}
{\em Let $M$ be a compact Riemannian surface with non-empty strictly convex boundary. Suppose that $V$ is an integral $1$-varifold in $M$ which is stationary with free boundary, and $\mathbb{Z}_2$-almost minimizing in small annuli with free boundary. Then $V$ is a free boundary geodesic network and such that $\theta^1(V,p)$ is a positive integer for all $p\in spt(||V||)\cap int(M)$. Moreover, if $p \in int(M)$ and $\theta^1(V,p)=2$, then, $p$ is a crossing point of two geodesic segments with multiplicity one.}

The $\mathbb{Z}_2$-almost minimizing condition is a fundamental variational property of objects that arise in min-max theory. In the present setting, the notion that is considered is that of \cite{Don}. In order to state this definition, recall that $\mathcal{Z}_{1,rel}(M, \partial M; \mathbb{Z}_2)$ represents the space of relative $1$-dimensional flat cycles modulo $2$. The relevant topologies in this space are those induced by the flat metric, the ${\bf F}$-metric, and the mass. See Section 2 of \cite{Don}, Section 2 of \cite{LioMarNev}, and Section 2 of \cite{marques-neves-infinitely}, for these and other relevant definitions. 

\subsection{Definition} Let $U$ be a relatively open subset of $M$. The $1$-varifold $V$ is $\mathbb{Z}_2$-almost minimizing in $U$ if for every $\varepsilon>0$, there exists $\delta>0$ and a relative cycle $S \in \mathcal{Z}_{1,rel}(M, \partial M; \mathbb{Z}_2)$ such that ${\bf F}(|S|,V)< \varepsilon$ and with the property that for every finite sequence $\{S_i\} \subset \mathcal{Z}_{1,rel}(M, \partial M; \mathbb{Z}_2)$, $0\leq i\leq m$, with $S= S_0$ and satisfying
$$
spt(S-S_i)\subset U,\ \ \mathcal{F}(S_i - S_{i-1})\leq \delta, \text{ and } {\bf M}(S_i)\leq {\bf M}(S) + \delta,
$$
for all $1 \leq i \leq m$, then it must hold that ${\bf M}(S_m) \geq {\bf M}(S) - \varepsilon$.

Next, let us recall the basic min-max set up, following \cite{Don}. Fix a cubical subcomplex $X$ of the $m$-dimensional cube $I^m = [0,1]^m$. A continuous map in the flat topology $\Phi : X \rightarrow \mathcal{Z}_{1,rel}(M, \partial M; \mathbb{Z}_2)$ is a sweepout if $\Phi^{\ast}(\overline{\lambda})\neq 0$ in $H^1(X; \mathbb{Z}_2)$, where $\overline{\lambda}$ is the generator of $H^1(\mathcal{Z}_{1,rel}(M, \partial M; \mathbb{Z}_2); \mathbb{Z}_2) \simeq \mathbb{Z}_2$ (by Almgren's isomorphism). Only sweepouts with no concentration of mass are considered in the definition of the min-max width. More precisely, let $\mathcal{P}$ be the set of sweepouts $\Phi$ of $M$ such that $\lim_{r\rightarrow 0} m(\Phi,r) = 0$, where
$$
m(\Phi,r) = \sup \{||\Phi(x)||(B(p,r)\setminus (\partial M)) : x \in dmn(\Phi) \text{ and } p\in M \},
$$
and $B(p,r)$ is the geodesic ball of $M$ centered at $p$ of radius $r$. Observe that the domain $dmn(\Phi)$ of admissible sweepouts $\Phi \in \mathcal{P}$ is not fixed.

\subsection{Definition} The (first) min-max width is the number
$$
\omega(M) = \inf_{\Phi \in \mathcal{P}} \left( \sup\{{\bf M}(\Phi(x)) : x \in dmn(\Phi)\}\right).
$$

Corollary $4$ of \cite{Don} implies that $\omega(M)$ can be arbitrarily approximated by the mass of $1$-varifolds which are free boundary stationary and $\mathbb{Z}_2$-almost minimizing in small annuli. Each varifold in the approximating sequence arises as solution of a standard homotopy min-max problem, see Theorem $3$ of \cite{Don}. The refinement of the combinatorial argument in the min-max theory of Almgren and Pitts obtained by Colding and De Lellis gives that the $\mathbb{Z}_2$-almost minimizing property holds either in small balls centered at all points of $M$, or in the complement of sufficiently small balls centered at a single point $p \in M$. This is explained in the proof of Proposition $5.1$ of \cite{ColDeL}. Combining this fact with Theorem \ref{structure-fb}, one has the following result.

\subsection{Theorem}\label{thm-single-jct}
{\em Let $M$ be a compact Riemannian surface with non-empty strictly convex boundary. Then, there exists a free boundary geodesic network $V$ with at most one junction such that $\omega(M) = ||V||(M)$. Moreover, if $V$ has an interior junction at $p\in M$, then the density $\theta^1(V,p)$ of $V$ at $p$ is a positive integer. If, in addition, it happens that $\theta^1(V,p)=2$, then $p$ is a crossing point of two geodesic segments with multiplicity one.}

\begin{proof}
Let $\{V_i\}$, $i\in \mathbb{N}$, be a sequence of free boundary geodesic networks in $M$, each with at most one junction, with ${\bf M}(V_i) \rightarrow \omega(M)$. Next, it will be explained that a subsequence of $\{V_i\}$ converges to a free boundary stationary geodesic network $V$ with at most one junction, such that ${\bf M}(V) = \omega(M)$. 

In order to achieve the goal, one proves that the density $\theta^1(V_i, p_i)$ at the junction must be uniformly bounded. This implies the previous claim, since it suffices to consider a subsequence such that the junctions converge to a point $p \in M$, and the lists of unit velocity vectors of geodesic segments of $V_i$ in $T_{p_i}M$ converge to a list of unit vectors in $T_p M$. Indeed, each $V_i$ is determined by the point $p_i$ and a finite list of vectors in $T_{p_i}M$, possibly with repetitions, which are the velocity vectors of the geodesic segments of $V_i$ emanating from the singular point $p_i$. Assuming that $\theta^1(V_i, p_i)$ are uniformly bounded, there exists a subsequence for which the densities are all equal to a fixed quantity. Up to a new subsequence, the vectors in these lists of a fixed number of vectors converge to a list with the same number of vectors in $T_pM$. Then, $V_i$ converges to the union of the geodesic segments in $M$ in the directions of those limit vectors pointing inside $M$.

It remains to show that the densities are bounded. Suppose, after passing to a subsequence, that the junctions converge to a point $p\in M$. If the limit point is in the interior of $M$, then, having unbounded densities would immediately imply that ${\bf M}(V_i) \rightarrow \infty$. Which contradicts ${\bf M}(V_i) \rightarrow \omega(M)$.

From now on, suppose that $p_i \rightarrow p \in \partial M$, with $\theta^1(V_i, p_i)\rightarrow \infty$. Let $\varepsilon_i = d(p_i, \partial M)$, where $d$ represents distance in $M$. The boundary $\partial M$ has positive geodesic curvature with respect to the inward pointing unit vector $\nu$ normal to $\partial M$. Let $\{\partial M_{t}\}$, $t \in [0, \varepsilon]$, be the foliation of a tubular neighborhood of $\partial M$ by convex curves, where $M_{\varepsilon} = \{x \in M : d(x, \partial M)\geq \varepsilon\}$. The opening angle of a vector $v \in T_{p_i}M$ is the shortest angle determined by this vector and the inward pointing unit vector $\nu_{i}$ normal to $\partial M_{\varepsilon_i}$.

Let $\{q_1, q_2\} = (\partial B(p, \varepsilon))\cap(\partial M)$, and let $\tau(p, q_j)$ be the shortest geodesic from $p$ to $q_j$, $j=1,2$. The maximum opening angle $\theta_0$ of the velocities of $\tau(p, q_j)$ at $p$ with respect to $\nu(p)$ are strictly smaller than $\pi/2$. Any geodesic $\sigma$ from $p = \sigma(0)$ with velocity $\sigma^{\prime}(0)$ of opening angle bigger than $\theta_0$ and smaller than $\pi/2$ is constrained to one of two narrow regions bounded by $\tau(p, q_j)$. Moreover, for small $\varepsilon$, the final point of $\sigma$ belongs to one of the arcs bounded by $p$ and one of the points $q_j$ making an angle with the boundary of $M$ which is far from $\pi/2$. By continuity, for large $i$, the geodesics leaving $p_i$ with opening angle bigger than $\theta_0$ end at a boundary point contained in a small neighborhood of the arc of $\partial M$ from $q_1$ to $q_2$ inside $\partial M$. The only such geodesic which meets $\partial M$ orthogonally is that with initial velocity $-\nu_{i}(p_i)$.

In particular, since each $V_i$ has a single junction, the velocity of any geodesic segment of $V_i$ at $T_{p_i}M$ is either equal to $-\nu_{i}(p_i)$, or has opening angle with respect to $\nu_i$ smaller than $\theta_0$. It follows from $\theta_0 < \pi/2$, $\theta^1(V_i, p_i)\rightarrow \infty$, and the stationarity condition at $p_i$, that the number of velocity vectors of geodesic segments of $V_i$ at $p_i$ with opening angle smaller than $\theta_0$, counted with repetitions, also grows to infinity as $i\rightarrow \infty$. Which implies that ${\bf M}(V_i) \rightarrow \infty$, since the lengths of all such geodesics are uniformly bounded from below away from zero.

The claims made in the theorem about the density of the limit $V$ at a hypothetical interior junction follow from the fact that density can not drop in the limit. The regularity at interior junctions of density $2$ is a general fact which follows from the stationarity of $V$ only, see \cite{Don}.
\end{proof}

\section{Applying the process of Birkhoff to sweep regions out}\label{Sec-Birkhoff}

Let $M$ be a compact Riemannian surface with strictly convex boundary. Let us fix a positive integer $L$ and $\varepsilon >0$. The unit interval will be denoted by $I = [0,1]$. Let us use $\Lambda$ to denote the set of broken geodesics $\sigma : I \rightarrow M$ parametrized with constant speed, with exactly $L-1$ break points (possibly with unnecessary breaks), $\sigma(0)$, $\sigma(1) \in \partial M$, each geodesic segment of $\sigma$ having length at most $\varepsilon$,  and with Lipschitz constant bounded by $L$.

Consider the subset $G$ of the immersed geodesics with free boundary. There exists an $\varepsilon >0$, depending on the geometry of $M$, for which a curve shortening map $D : \Lambda \rightarrow \Lambda$ with the following properties can be constructed:
\begin{enumerate}
    \item[(a)] $D$ is continuous with respect to the $W^{1,2}$-norm;
    \item[(b)] $D(\sigma)$ is homotopic to $\sigma$;
    \item[(c)] $L(D(\sigma))\leq L(\sigma)$, for all $\sigma \in \Lambda$;
    \item[(d)] If $D(\sigma) = \sigma$, then $\sigma \in G$;
    \item[(e)] $L(D(\sigma))< L(\sigma)$, for all $\sigma \in \Lambda \setminus G$.
\end{enumerate}

The details of this construction and further properties are explained in \cite{Zhou-FB-Geod}. This shortening map is an adaptation for the boundary case of the Birkhoff curve shortening process. The above configuration is included in those of \cite{Zhou-FB-Geod} by embedding $M$ isometrically in a closed Riemannian surface and choosing $N = \partial M$, in the notation of that article. The convexity implies that the deformations of broken geodesic in $M$ do not leave this domain. 

The construction can be briefly explained as performing two geodesic approximations using a set of $2L$ evenly spaced points $\{\sigma(x_i)\}$, $1\leq i \leq 2L$, in $\sigma$. The first step consists of replacing the geodesic segments $\sigma|_{[0,x_2]}$ and $\sigma|_{[x_{2L-2},x_{2L}]}$ that arrive at the boundary points $\sigma(0)$ and $\sigma(1)$ by shortest geodesic segments connecting the points $\sigma(x_2)$ and $\sigma(x_{2L-2})$ to $\partial M$, respectively. And then, replacing $\sigma$ on each interval $[x_{2i}, x_{2i+2}]$ by the shortest geodesic segment with the same endpoints. This finishes the first geodesic approximation. After a reparametrization to constant speed, one performs a second geodesic approximation by conecting the midpoints of the geodesic segments of the curve obtained after the first approximation.

Let $\Omega$ be a compact subset of $M$ such that its interior is diffeomorphic to an open $2$-disk. Assume that $\partial \Omega$ is the union of a broken geodesic and a segment of $\partial M$ , possibly empty or degenerated to a point. More precisely, the boundary of $\Omega$ is the image of a closed curve $\alpha : [0,2] \rightarrow M$ such that $\alpha|_{[0,1]}$ is a broken geodesic, and $\alpha|_{[1,2]}$ parametrizes a portion of $\partial M$. The broken geodesic is not necessarily simple; there might be $t_1 \neq t_2$ such that $\alpha$ touches itself at $\alpha(t_1) = \alpha(t_2)$ on the exterior region of $\Omega$, without crossing. This includes the case in which a geodesic segment of $\alpha$ can be approximated by points in $\Omega$ from both sides, i.e., two segments of $\alpha$ parametrized by disjoint intervals in $[0,1]$ have exactly the same image. Even though $\alpha([0,1])$ is a broken geodesic, it may touch the boundary of $M$ finitely many times. 

Suppose, in addition, that the inner angles of $\Omega$ are all smaller than $\pi$, and that the angles between $\alpha|_{[0,1]}$ and $\alpha|_{[1,2]}$ measure at most $\frac{\pi}{2}$.

The following lemma is the main consequence of the Birkhoff process that will be applied in Section \ref{Sec-W-type}. It is the natural extension of Lemma 1.2 in Part II of \cite{CalCao} to the present setting. See also Lemma 2.2 in \cite{Cro}. The result follows by considering the iterated sequence of deformations $\{\alpha_i\}$ defined as $\alpha_1 = \alpha|_{[0,1]}$ and $\alpha_{i} = D(\alpha_{i-1})$, for $i\geq 2$. Once the construction adapted to the boundary is performed, only minor changes to the proof need to be pointed out. This is explained in the paragraphs after the statement. 

\subsection{Lemma}\label{varrer-convexos}
{\em Let $M$, $\Omega$ and $\alpha : [0,2]\rightarrow M$ be as above. Suppose that $\alpha|_{[0,1]}$ is not a free boundary geodesic (without break points), and that $\alpha|_{[1,2]}$ is not empty or a constant curve. Then, there exists either a simple free boundary geodesic $\gamma$ of length $L(\gamma) < L(\alpha|_{[0,1]})$ and contained in the relative interior of $\Omega$, or a homotopy $\{\sigma_s\}$, $s\in I$, which satisfies:
\begin{enumerate}
    \item[(1)] $\sigma_s$ is a broken geodesic with $\sigma_s(0)$,$\sigma_s(1) \in (\partial M)\cap \overline{\Omega}$, for all $s\in I$.
    \item[(2)] $\sigma_0 = \alpha|_{[0,1]}$ and $\sigma_1$ is a constant curve.
    \item[(3)] $L(\sigma_s) \leq L(\alpha|_{[0,1]})$, for all $s\in I$.
    \item[(4)] $\sigma_s$ may have self-intersections which coincide with those in $\alpha([0,1])$. In particular, these are exterior self-intersections. 
    \item[(5)] $\sigma_s$ has convex angles with respect to $\Omega_s = \Omega \setminus \{\sigma_t : 0\leq t \leq s\}$.
    \item[(6)]
    $\sigma_s\rightarrow \sigma_{s_0}$ continuously as graphs, for all $s_0 \in I$.
    \item[(7)]
    $\sigma_s\rightarrow \sigma_{s_0}$ smoothly as graphs away from break points, for all $s_0 \in I$.
\end{enumerate}}

Besides the adaptation due to the boundary, Lemma \ref{varrer-convexos} differs from the versions found in \cite{CalCao} and \cite{Cro} because of possible exterior self-intersections.

The convexity property described in Definition 1.0 of Part II of \cite{CalCao} can be replaced by the following notion: there exists an $\varepsilon >0$ such that for all $x$ and $y$ in $\alpha([0,1])$ with $d_{\Omega}(x,y)< \varepsilon$, the minimizing geodesic $\tau$ from $x$ to $y$ satisfies $\tau \subset \overline{\Omega}$, where $d_{\Omega}(x,y)$ is the infimum of lengths of piecewise smooth curves $c: [0,1] \rightarrow \overline{\Omega}$ from $x$ to $y$, and such that $c(t) \in \Omega$, for all $0<t<1$. 

The non-constant portions of $\alpha$ parametrized by a subinterval $[a,b]\subset [0,1]$ and connecting a point of self-intersection with itself, $\alpha(a) = \alpha(b)$, can not be contained in a ball of radius smaller than the injectivity radius of $M$ and centered at $\alpha(a)$, unless $\Omega$ is contained in such a ball. Indeed, if the loop $\alpha([a,b])$ is contained in a ball of sufficiently small radius, it would divide the ball in two components and would be convex with respect to one of the sides. By the convexity properties of small geodesic balls, this loop can not be convex to its exterior. Therefore, $\Omega$ must be contained in that ball. However, since $\alpha([1,2])$ is not constant, this case is ruled out. 

In particular, for each $x\in \alpha([0,1])$, the longest open connected segment $\alpha_x$ of $\alpha|_{[0,1]}$ that contains $x$ and $d(x,y)< inj(M)$, for all $y \in \alpha_x$, does not have self-intersections. Next, the convexity assumption mentioned above gives us that the minimizing geodesic $\tau$ from $x$ to $y$ satisfies $\tau \subset \overline{\Omega}$, and either $\tau\cap \alpha([0,1]) = \{x,y\}$ or $\tau \subset \alpha([0,1])$. This is the analog of item (1) of Lemma 2.1 of \cite{Cro} adapted to the present setting.  

Moreover, the assumption on the angles of $\alpha$ at $\alpha(0)$ and $\alpha(1)$ measuring at most $\frac{\pi}{2}$ implies that: if $x \in \alpha([0,1])$ and $d(x,\alpha([1,2]))< inj(M)$, then the shortest geodesic $\tau$ from $x$ to $\partial M$ satisfies $\tau \subset \overline{\Omega}$, and either $\tau \cap \alpha([0,1]) = \{x\}$ or $\tau \subset \alpha([0,1])$. The remaining of the proof of Lemma 2.2 of \cite{Cro}, implies that self-intersections along the Birkhoff curve shortening process that is used here only occur at points of self-intersection of $\alpha([0,1])$. The convexity assumption also follows as in \cite{Cro}. The iteration applied in Lemma 1.2 of Part II of \cite{CalCao} ends the argument. It should be noted that the number of breaks in the Birkhoff process depends on $L(\alpha([0,1]))$, as in \cite{CalCao} and \cite{Cro}.

\subsubsection*{Remark:} If the assumption that $\alpha|_{[1,2]}$ is not empty or a constant curve is replaced by its opposite, then $\alpha([0,1])$ is closed and the usual Birkhoff argument implies that either $\Omega$ has a closed geodesic in its interior, or there is a homotopy as above by closed broken geodesics.

\section{Geodesic networks: properties and induced sweepouts}\label{Sec-W-type}

Let $g$ be a Riemannian metric in a surface $M$ that is diffeomorphic to the closed $2$-dimensional disk. Assume that the sectional curvature of $g$ is nonnegative, $K\geq 0$, and that its boundary is strictly convex. 

Let $V$ be a free boundary geodesic network, as in Definition \ref{FBGN}. In addition, assume that $V$ satisfies the following condition:
\begin{enumerate}
    \item[($D_i$)] $\theta^1(V, p)$ is a positive integer, for all $p \in spt(||V||)\cap int(M)$.
\end{enumerate}
For fixed $V$, let $W$ be a varifold obtained from $V$ after the deletion of some, possibly none, geodesic segments. Assume that $W$ satisfies:
\begin{enumerate}
    \item[($w_1$)] $\theta^1(W, p)$ is a positive integer, for all $p \in spt(||V||)\cap int(M)$.
    \item[($w_2$)] the angles between two consecutive geodesic segments of $W$ at an interior junction are strictly smaller than $\pi$.
    \item[($w_3$)] the smallest angle between any vector in $T_p (\partial M)$ and the velocity vector of a geodesic segments of $W$ emanating from a boundary junction $p \in \partial M$ is always smaller than or equal to $\frac{\pi}{2}$.
\end{enumerate}

One can check that $V$ satisfies ($w_1$), ($w_2$), and ($w_3$), for instance. The first result in this section is a refined structure characterization for varifolds such as $V$ and $W$, which is a consequence of the non-negative curvature and convexity assumptions. Consider the following definition.

\subsection{Definition}
Let $\Omega$ be a compact subset of $M$ such that its interior is diffeomorphic to an open $2$-disk. $\Omega$ is said to satisfy the $\ast$-property if 
\begin{enumerate}
    \item[(a)] $\partial \Omega = \alpha([0,2])$, where $\alpha : [0,2] \rightarrow M$ is a closed curve such that $\alpha|_{[0,1]}$ is a broken geodesic, and $\alpha|_{[1,2]}$ parametrizes a portion of $\partial M$, which is possibly empty or degenerated to a point. The portion $\alpha([0,1])$ might touch itself on the exterior region, as explained in Section \ref{Sec-Birkhoff}.
    
    \item[(b)] The inner angles of $\Omega$ are all smaller than $\pi$, and that the angles between $\alpha([0,1])$ and $\alpha([1,2])$, if they exist, measure at most $\frac{\pi}{2}$.
\end{enumerate}

\subsection{Lemma}\label{W-properties}
 {\em Let $V$ and $W$ be as above, and let $\{\Omega_i\}$, $1\leq i \leq q$, be the connected components of $M\setminus spt(||W||)$, where $spt(||W||)$ represents the support of the weight measure of $W$. Then, each $\Omega_i$ is diffeomorphic to an open $2$-disk, $spt(||W||)$ intersects $\partial M$ and is connected. Each $\Omega_i$ satisfies the $\ast$-property. Moreover, there exists a decomposition of $\{\Omega_i\}$ as the union of disjoint collections $\{\Omega_j\}_{j\in \mathcal{I}} \cup \{\Omega_k\}_{k\in \mathcal{J}}$ with $\mathcal{I} \cup \mathcal{J} = \{i : 1\leq i 
 \leq q\}$ and
 \begin{equation}\label{decomp}
 \sum_{j\in \mathcal{I}} [\partial \Omega_j] = \Gamma = \sum_{k\in \mathcal{J}} [\partial \Omega_k], 
 \end{equation}
 where $[\partial \Omega_j]$ denotes the induced relative flat chain modulo $2$, and $\Gamma$ denotes the sum of the geodesic segments of $W$ with odd multiplicity.}

\begin{proof}
Observe that hypothesis ($w_2$) implies that the turning angles $\theta_i$ of a connected component $\Omega$ of $M\setminus spt(||W||)$ are all positive. Moreover, the boundary of $\Omega$ is a piecewise smooth curve made of smooth curves $c_i$ which are either geodesic segments of $W$, or arcs of the convex boundary $\partial M$. These observations and the Gauss-Bonnet formula,
$$
\int_{\Omega} K \hspace{.1cm} d\mathcal{H}^2 + \sum_i \int_{c_i} k_g\hspace{.1cm} d\ell + \sum_{i} \theta_i = 2\pi \chi(\Omega), 
$$
imply that $\chi(\Omega) \geq 0$. Here $\chi(\Omega) = 2 - 2g - r$, where $g = genus(\Omega) =0$ and $r$ is the number of boundary components of $\partial \Omega$. The fact that $g=0$ follows from the fact that $\Omega \subset M$, and $M$ is diffeomorphic to a disk. Therefore, we conlude that $r=1$ or $2$. We claim that $r=1$. Suppose, by contradiction, that $r = 2$. Then the Gauss-Bonnet formula implies that the boundary of $\Omega$ must be equal to two disjoint closed geodesics in $int(M)$. On the other hand, the Gauss-Bonnet formula applied to the region whose boundary is made of a closed geodesic and $\partial M$ gives a contradiction. Then, $r=1$.

Note that the same argument shows that if $W$ as above is non-trivial, $W \neq 0$, then each connected component of $spt(||W||)$ intersects $\partial M$.

Let us prove, by a contradiction argument, that the support of $W$ is connected. Assuming that $spt(||W||)$ had multiple components, let $W_i$, $1\leq i \leq N$, be the sum of the geodesic segments (with the same respective multiplicities as in $W$) contained in the support of each of those different components of $spt(||W||)$. Then, there would exist a connected component $\Omega$ of $M\setminus spt(||W||)$ such that $(\partial \Omega) \setminus (\partial M)$ would be disconnected. Indeed, a connected component $\Omega^{\ast}$ of $M\setminus spt(||W_1||)$ would contain the support of some other $W_i$, $i\geq 2$. Let $\Omega$ be a connected component of $\Omega^{\ast}\setminus spt(||W||)$ whose boundary touches $spt(||W_1||)$. This choice of $\Omega$ has the desired property, since its boundary intersects the supports of $W_1$ and some other $W_i$, which we assume, up to changing the indices, is $W_2$ . Note that $\Omega$ has at least four corners, two at the intersections of $\partial M$ with $spt(||W_1||)$, and two other at $spt(||W_2||)\cap (\partial M)$. Hypothesis ($w_3$) implies that the corresponding turning angles measure at least $\frac{\pi}{2}$. In addition, $\partial \Omega$ contains a non-degenerate segment of $\partial M$. Since $\sum_j \theta_j \geq 4\cdot \frac{\pi}{2} = 2\pi$, $\sum_i \int_{c_i} k_g\hspace{.1cm} d\ell > 0$, and $K\geq 0$, the Gauss-Bonnet formula applied to $\Omega$ gives us that $2\pi \chi(\Omega) > 2\pi$, which contradicts the fact that the interior of $\Omega$ is diffeomorphic to a disk.

In what follows, the decomposition $\mathcal{I} \cup \mathcal{J} = \{i : 1\leq i \leq q\}$ is explained. Fix $1\in \mathcal{I}$. For each $2\leq j \leq q$, let $\sigma$ be a path in $M$ joining a point $p_1 \in \Omega_1$ and $p_j \in \Omega_j$ which does not cross geodesic segments of $W$ at a singular point. The remainder modulo $2$ of the number of crossing points of $\sigma$ with $spt(||W||)$ (the multiplicity of the geodesic segments is considered, e.g., if $\sigma$ crosses a geodesic with multiplicity $2$, then two intersections are considered) is independent of the path. Indeed, this intersection number is homotopy invariant whenever the homotopy avoids the singular points. In addition, hypothesis ($w_1$) implies the invariance through the singular point. Indeed, before crossing the singular point $p$, we can perturb $spt(||W||)$ around a small neighborhood of this point to a union smooth curves crossing at $p$ without changing the intersection with $\sigma$. This is possible because of ($w_1$). Next, apply the homotopy starting at $\sigma$, and deform the union of crossing curves back to $W$ without changing the intersections with the final image of $\sigma$. It is now clear that the remainder modulo $2$ of the intersection number does not change through the homotopy. We define the decomposition putting $j\in \mathcal{I}$ if that remainder is zero, and $j \in \mathcal{J}$ if the remainder is one.

It follows from the construction that two $\Omega_j$ that share a common edge of $W$ with odd multiplicity belong to different sets of the decomposition, while components of $M\setminus spt(||W||)$ sharing an edge with even multiplicity belong to the same set. In the second case, this common edge does not appear in one of the sums of (\ref{decomp}) and is cancelled in the other.
\end{proof}

The next result consists of an application of Lemma \ref{varrer-convexos} to the broken geodesics that make part of each $\partial \Omega_i$ from the previous lemma.

\subsection{Proposition}\label{varrer-regioes}
{\em Let $V$, $W$ and $\Omega_i$ be as in the statemet of Lemma \ref{W-properties}. Let us write $W = \sum_{j} m_j \gamma_{j}$ as the sum of its geodesic segments $\gamma_j$ with multiplicities $m_j$ which meet at junctions or at the boundary of $M$. There exists a sweepout of $M$ by curves of lengths bounded by the maximum between
$$
\sum_{j\in \mathcal{O}} L(\gamma_j) + 2\cdot \sum_{j\in E_{\mathcal{I}}} L (\gamma_j) \text{ and } \sum_{j\in \mathcal{O}} L(\gamma_j) + 2 \cdot \sum_{j\in E_{\mathcal{J}}} L (\gamma_j),
$$
where $\mathcal{O}$ represents the subset of indices for which $m_j$ is odd, $E_{\mathcal{I}}$ the subset of indices with even $m_j$ and such that $\gamma_j$ is contained in the closure of some $\Omega_i$ in $\mathcal{I}$, and $E_{\mathcal{J}}$ the set of indices which are not in $\mathcal{O} \cup E_{\mathcal{I}}$. 

If $W$ contains a free boundary geodesic $\gamma$, then $M$ can be swept out by curves of lengths bounded by $L(\gamma)$. If, additionally, the curve $\gamma$ has one self-intersection in the interior of $M$, then one can sweep $M$ out by curves whose supremum of lengths is strictly lower than $L(\gamma)$.}

\begin{proof} Let us assume first that the support of $W$ contains a free boundary geodesic $\gamma : [0,L] \rightarrow M$, parametrized by arc length, which connects two points of $\partial M$, and without self-intersections. Let $H: [0,L] \times (-\varepsilon, \varepsilon) \rightarrow M$ be a smooth variation of $\gamma = H(\cdot, 0)$ such that $H(0,t), H(L,t) \in \partial M$, for all $-\varepsilon <t < \varepsilon$. The variational vector field $V = \frac{\partial H}{\partial t}(\cdot,0)$ satisfies $V(0) \in T_{\gamma(0)}(\partial M)$ and $V(L) \in T_{\gamma(L)}(\partial M)$. Choose $V=N$, the unit vector field perpendicular to $\gamma$. This is possible because of the free boundary condition. Since $\gamma$ is a free boundary geodesic, the first variation of the lengths $L(H(\cdot, t))$ at $t=0$ vanishes. Letting $A_{\partial M}$ represent the second fundamental form of $\partial M$, and $K(p)$ be the sectional curvature of $M$ at point $p$, the second variation formula gives:
$$
\frac{d^2}{dt^2}L(H(\cdot, t))\big|_{t=0} = - A_{\partial M}(V(0), V(0)) - A_{\partial M}(V(L), V(L)) - \int_0^L K(\gamma(s)) ds.
$$
The convexity and $K\geq 0$ assumptions imply that this deformation has negative second derivative. Therefore, reducing $\varepsilon$ if necessary, the length strictly decreases in $[0,\varepsilon)$, and increases in $(-\varepsilon, 0]$.

Homotope the curve $H(\cdot,\varepsilon)$ to a broken geodesic approximation $\alpha_0$ with a sufficiently large number of breaks depending on the length $L(H(\cdot,\varepsilon))$. Observe that $L(\alpha_0)< L(H(\cdot,\varepsilon)) < L(\gamma)$ and $\alpha_0$ is contained in a connected component of $M \setminus \gamma$. Apply the Birkhoff's curve shortening process $D$ to obtain $\alpha_i = D(\alpha_{i-1})$, for $i\geq 1$. These curves all have extremities on $\partial M$, and they converge either to a free boundary geodesic, or to a constant curve. Since $\gamma$ is a free boundary geodesic, it serves as a barrier to the curves $\alpha_i$. This implies that $\alpha_i$ cannot converge to a free boundary geodesic, otherwise, this geodesic would be disjoint from $\gamma$ and their union would contradict the connectedness obtained in Lemma \ref{W-properties}. Hence, the $\alpha_i$ converge to a constant curve, and the iterated deformations originate a flat continuous path with no concentration of mass from $\alpha_0 = H(\cdot, \varepsilon)$ to a constant curve. This path is slightly different from those obtained in Lemma \ref{varrer-convexos} because the convexity assumption may fail at $\alpha_0$. This is not the case if the sectional curvature is strictly positive along $\gamma$. It turns out that the convergence properties (6) and (7) still hold, and one can combine a similar path from $H(\cdot, -\varepsilon)$ and interpolate to obtain a sweepout by curves of length at most the length of $\gamma$. Further details about the interpolation are given at the end of this proof. 

Assume now that the free boundary geodesic $\gamma$ has one self-intersection in $int(M)$, which is a single self-intersection forming a loop inside $M$. There are three geodesic segments emanating from this point of intersection only. One is the loop, and the other two intersect the boundary perpendicularly. The argument here is exactly as in the proof of Claim 1 in Theorem D, Part II, of \cite{CalCao}. As in the case explained above, one starts by using the second variation to obtain a shorter curve with the same shape of $\gamma$. Then, performing the exact same cut-and-paste constructions as in \cite{CalCao}, one obtains the desired sweepout by curves of lengths bounded by a number strictly below $L(\gamma)$. The only change here is that, in the present case, one uses Lemma \ref{varrer-convexos} instead of the corresponding result in the closed case. There are no changes in the cut-and-paste argument because it is performed in a small neighborhood of a point in $int(M)$. A similar analysis is explained in details in Proposition \ref{V-figure}. Next, we analyze the general case considered in the statement.

The properties of $\Omega_i$ obtained in Lemma \ref{W-properties} imply that there exists a parametrization $\alpha_i : [0,2] \rightarrow M$ of its boundary in such a way that $M$, $\Omega_i$, and $\alpha_i$ are either exactly as in the statement of Lemma \ref{varrer-convexos}, or $\alpha_i([0,1])$ is a closed broken geodesic and the closed case of the Birkhoff process can be applied. It follows from Lemma \ref{varrer-convexos}, or the equivalent result explained in the remark after that lemma, that there exist homotopies sweeping out each $\Omega_i$ by curves of length bounded by $L(\alpha_i([0,1]))$. Indeed, the only other option would be the existence of a free boundary (or closed) simple geodesic disjoint from $W$ and with length strictly less than $L(\alpha_i([0,1]))$. The assumptions $K\geq 0$ and convex boundary imply that Lemma \ref{W-properties} could be applied to the sum of $W$ with that simple geodesic. This would give a contradiction with the fact that a varifold $W$ in the statement of Lemma \ref{W-properties} must be connected.

Observe that geodesic segments of $W$ which bound $\Omega_i$ on both sides, are parametrized twice by $\alpha_i|_{[0,1]}$. The intersection number argument used in the proof of Lemma \ref{W-properties} implies that the multiplicity of each such geodesic segment must be even. All the other geodesic segments of $W$ in $\alpha_i$ are parametrized exactly once by $\alpha_i$. The geodesic segments $\gamma_j$ of $W$ with $j \in \mathcal{O}$ are parametrized exactly once by a curve $\alpha_i$ with $i \in \mathcal{I}$, and once by a curve $\alpha_j$ with $j\in \mathcal{J}$, where $\mathcal{I}$ and $\mathcal{J}$ are the sets of indices associated to the decomposition of the sets $\{\Omega_i\}$ made in the statement of Lemma \ref{W-properties}.

The homotopies $\{\sigma^i_s\}$, $s\in [0,1]$, associated to all $\Omega_i$ with $i \in \mathcal{I}$ are combined to form a continuous path, with respect to the topology of flat chains modulo $2$ relative to $\partial M$. This combination is performed by adding, for each fixed $s\in [0,1]$, the curves $\sigma^i_s$, for $i \in \mathcal{I}$, and it deforms the set $\Gamma$ given in the statement of Lemma \ref{W-properties} to a finite set in $\cup_{i\in \mathcal{I}} \Omega_i$. The continuity follows from property (6) listed in the statement of Lemma \ref{varrer-convexos}.

Performing the same construction with the homotopies supported in the $\Omega_j$, $j\in \mathcal{J}$, one obtains a flat continuous sweepout $\Phi : \mathbb{S}^1 \rightarrow \mathcal{Z}_{1}(M, \partial M;\mathbb{Z}_2)$, which satisfies the no concentration of mass hypothesis for interpolation: 
$$
\lim_{r\rightarrow 0} \sup \{||\Phi(t)||(B_r(x)) : t \in \mathbb{S}^1 \text{ and } x \in M\} = 0.
$$
This follows from the stronger convergence listed as property (7) of Lemma \ref{varrer-convexos}. Moreover, it follows from property (3) of Lemma \ref{varrer-convexos}, and the observations made above in this proof about the number of times that each segment is parametrized by the $\alpha_i$, that the mass of each $\Phi(t)$ is bounded from above by the quantity described in the statement of this proposition.

Finally, the continuous-to-discrete interpolation theorem, Theorem 2.10 of \cite{LioMarNev} and references therein, gives a discrete sweepout with the properties that were listed in the statement of the proposition.
\end{proof}

\section{Geodesic networks with a single junction}\label{sec-single-junction}

In this section, the attention is focused on special varifolds $V$ and $W$ that satisfy the properties stated in Section \ref{Sec-W-type}. Suppose that $V$ is not a single free boundary geodesic without self-intersections, or with a single crossing point in the interior of $M$. In addition, it is assumed that the multiplicity of each geodesic segment of $V$ is either $1$ or $2$, and that $V$ does not contain a simple free boundary geodesic. Otherwise, Proposition \ref{varrer-regioes} would imply that $M$ could be swept out by flat chains modulo $2$ of masses bounded by a number strictly lower than the mass of $V$. More importantly, assume that $V$ has a unique singular point. More precisely, consider the following definition.

\subsection{Definition}\label{admissible}
Let $V$ be a free boundary geodesic network satisfying ($D_i$). One says that $V$ has an \textit{admissible singular point} at $p\in M$ if $p$ is a singular point, i.e., a point $p$ around which $V$ is not a single geodesic (orthogonal to the boundary of $M$ at $p$ in the case that this point belongs to $\partial M$), which belongs to one of the following classes:
    \begin{enumerate}
        \item[($J_i$)] $p$ is an interior junction at which $\theta^1(V,p) \geq 3$;
        \item[($J_b$)] $p$ is a boundary junction at which $\theta^1(V,p)\geq 1.5$;
        \item[($J_{\ell}$)] $p \in \partial M$, $\theta^1(V,p)= 1$, and $V$ is not a geodesic loop.
    \end{enumerate}

Observe that the case in which $p \in int(M)$ and $\theta^1(V,p) = 2$ is not included in the definition, since Theorem \ref{structure-fb} implies that it would be the crossing point of two local geodesics. But, in this case, $V$ would necessarily contain a free boundary geodesic with at most one self-intersection.  

Assume that each geodesic segment of $V$ starts and ends at points which are either of type ($J_i$), or ($J_b$), or ($J_{\ell}$), or a regular point at the boundary of $M$. In the last case, the geodesic segment must be perpendicular to $\partial M$. This can be assumed simply by combining geodesic segments that meet at interior regular points of $V$, possibly with their respective multiplicities.

The uniqueness assumption narrows the analysis down to the situation at which each geodesic segment of $V$ is of one of the following types:
\begin{enumerate}
    \item[(i)] a geodesic segment that meets the boundary of $M$ orthogonally at one of its extremities and the singular point at the other;
    \item[(ii)] a geodesic loop starting and ending at the singular point. In this case, the inner angle of the loop measures strictly less than $\pi$.  
\end{enumerate}
A loop is a geodesic segment starting and ending at the same point of $M$, which will be referred to as the vertex of the loop. Assume that a varifold $W$ obtained from $V$, as described in the beginning of Section \ref{Sec-W-type}, has a loop. The inner angle of this loop is the angle determined by the two velocity vectors of the geodesic segment at the vertex, the initial and final velocities, with respect to the region of $M$ whose boundary is this loop. The velocity vectors at the vertex are both considered as pointing inside the loop. In the case that the vertex of the loop is at $\partial M$, the inner angle is immediately lower than $\pi$ because of the convexity assumption. Indeed, this angle can not be $\pi$ because $\partial M$ is not a geodesic. In the case of a vertex $p\in int(M)$, this angle property is implied by the Gauss-Bonnet formula applied to the region whose boundary is the union of the loop and the boundary of $M$.

Summarizing, $V$ is a free boundary geodesic network satisfying:
\begin{itemize}
    \item $V$ does not contain a free boundary geodesic with at most one self-intersection ($\theta^1(V,p)\geq 3$, if there exists a point of self-intersection);
    \item $V$ has a unique singular point $p$, which is admissible;
    \item if $p\in int(M)$, then $\theta^1(V, p)$ is a positive integer greater than two;
    \item each geodesic segment is either a type-(i) segment or a type-(ii) loop;
    \item the multiplicity of each geodesic segment is either $1$ or $2$.
\end{itemize}
In addition, $W$ is a varifold obtained from $V$ after the deletion of some geodesic segments in such a way that properties ($w_1$), ($w_2$), and ($w_3$) hold. The deletion could simply mean to drop the multiplicity of a given segment.

If the singular point $p$ is in the interior of $M$, two geodesic segments $\gamma_1$ and $\gamma_2$ of $W$ determine two angles at $p$. These are the angles between one velocity vector $v_1$ of $\gamma_1$, and a $v_2$ of $\gamma_2$ both emanating from $p$, such that no other velocity vector at $p$ from either $\gamma_i$, $i=1,2$, is inside these angles. In other words, these are the angles determined by consecutive inner velocity vectors of $\gamma_1$ and $\gamma_2$ which are not an inner angle of a loop. The vectors $v_1$ and $v_2$ belong to $T_pM$, which is considered with the inner product induced by the Riemannian metric. If both segments are of type-(i), for example, then exactly one of these angles measures strictly less than $\pi$. Indeed, since $V$ does not contain free boundary geodesics, it is impossible that both angles determined by $\gamma_1$ and $\gamma_2$ in this case measure $\pi$. In case one or both $\gamma_i$ are loops, then it is possible that the measure of the two angles determined by these segments are lower than $\pi$.

In the case that the singular point $p$ is at the boundary, only one angle determined by two geodesic segments is considered, namely, the angle which contains tangent vectors at $p$ that point inside $M$ only. In addition, the two angles between a given geodesic segment of $W$ and $\partial M$ are also considered. Since $V$ does not contain simple free boundary geodesics, then exactly one of the two angles determined by a type-(i) geodesic segment and the boundary measures less than $\frac{\pi}{2}$.

In what follows, Proposition \ref{varrer-regioes} is applied to show that one can sweep $M$ out by curves whose supremum of lengths is less than the total mass of $W$, unless additional restrictions on the multiplicities and angles are made.

\subsection{Proposition}\label{Prop-multiplicities-angles}
{\em Let $V$ and $W$ be as above, and assume that one of the following configurations hold.
\begin{enumerate}
\item[(a)] $W$ contains a geodesic loop inside another loop $\gamma$; letting $U\subset M$ be the connected region bounded by $\gamma$, it is assumed that $W$ contains another loop inside $\overline{U}$. Which can be $\gamma$ if it has multiplicity $2$ in $W$.

\item[(b)] The singular point $p$ is in $int(M)$, and $W$ contains two geodesic segments $\gamma_1$ and $\gamma_2$ which determine an angle $\theta$ measuring less than $\pi$. Assume that $W$ contains either a type-(ii) loop, or a pair of type-(i) segments whose velocities emanating from $p$ belong to $\theta$. The latter includes copies of the $\gamma_i$, $i=1,2$, if they are both of type-(i) and have multiplicity $2$, or a single type-(i) segment with multiplicity $2$ in $W$ with interior velocities at $p$ in the interior of the angle $\theta$.

\item[(c)] The singular point $p$ is in the boundary of $M$, and $W$ contains at least three geodesic  segments, counting multiplicities. It could be one segment with multiplicity one and another with multiplicity two.
\end{enumerate}
Then it is possible to sweep $M$ out by curves whose supremum of lengths is strictly lower than the total mass of $W$.}

\begin{proof} This is a consequence of Proposition \ref{varrer-regioes}. In cases (a) and (b), let $W^{\prime}$ be the varifold obtained from $W$ after deletion of the additional segments that are assumed to exist. In case (a), $W^{\prime}$ is obtained from $W$ by deleting the extra loop inside $\gamma$. If the original multiplicity of $\gamma$ in $W$ was two, then $W^{\prime}$ could be considered as simply lowering the multiplicity of $\gamma$ from two to one. In case (b), delete the type-(ii) loop or the pair of type-(i) segments.

In case (c), if there is a loop inside a loop, this case is reduced to (a). Otherwise, let $\gamma_1$, $\gamma_2$ and $\gamma_3$ be three of the segments that are assumed to exist. The velocity vectors of the three segments would be ordered in $T_pM$. More precisely, as one continuously rotates the vectors from one unit vector of $T_p(\partial M)$ to the other passing exactly once through each unit vector that points inside $M$, one would meet all the inner tangent vectors of $\gamma_1$ at $p$, then all vectors of $\gamma_2$, and, finally, the inner tangent vectors of $\gamma_3$. This ordering holds because if the vectors were alternating, then one segment would have two inner tangent vectors, then it would be a loop, and a tangent vector of another segment would be trapped inside this loop. The only possibility in this case would be a loop inside a loop. Therefore, the ordering of the velocities is proved. In this situation, let $W^{\prime}$ be obtained from $W$ after the deletion of the segment $\gamma_2$, which is not adjacent to $\partial M$.   

By assumption, $||W^{\prime}||(M) < ||W||(M)$. Moreover, it is easily checked that $W^{\prime}$ satisfies properties ($w_1$), ($w_2$), and ($w_3$). Indeed, the shortest angles between $W^{\prime}$ and $\partial M$ are exactly the same as the angles that this boundary makes with $W$. Only geodesic segments inside the loop were deleted in case (a). In (b), the deleted segments either intersect $\partial M$ at regular points only, or do not intersect $\partial M$ at all. In (c), only the middle segment was deleted.

Property ($w_1$) is kept because the density at the vertex is lowered exactly by the multiplicity of a loop in case (a), or by one in (b). The assumptions that the angles determined by the geodesic segments in (b) measure less that $\pi$ imply ($w_2$). In case (a), it is a consequence of the angle property listed in the characterization of type-(ii) loops. And ($w_2$) follows from the convexity of the boundary of $M$ in case (c). Therefore, Proposition \ref{varrer-regioes} implies that $M$ can be swept out by curves of lengths bounded by $||W^{\prime}||(M)$. 
\end{proof}

Next, the assumption that the unique singular point is admissible, see Definition \ref{admissible}, will be combined with Proposition \ref{Prop-multiplicities-angles}. The goal is to show that $V$ with a single admissible junction does not realize the min-max width.

\subsection{Corollary}\label{bdry}
{\em Let $V$ be as above with a singular point of type ($J_b$); a junction at the boundary such that $\theta^1(V,p) \geq 1.5$. Then, $M$ can be swept out by curves of lengths bounded by a number strictly lower than $||V||(M)$.}

\begin{proof}
In the case that $V$ contains at least three geodesic  segments, counting multiplicities, item (c) of Proposition \ref{Prop-multiplicities-angles} gives the desired property.

From now on, assume that $V$ has exactly two geodesic segments, counting multiplicities. Since $\theta^1(V,p)\geq 1.5$, the case of two geodesic segments which are not loops is not included in this analysis. 

Assume first that $V$ has one type-(i) segment $\gamma_1$, and one type-(ii) loop which will be denoted by $\gamma_2$. The curve $\gamma_1$ determines a unique angle $\theta_1$ measuring less than $\frac{\pi}{2}$ with a vector in $\partial M$. The stationarity of $V$ implies that the tangent vectors of $\gamma_2$ at $T_pM$ emanating from $p$ are not contained in $\theta_1$. Observe that $M\setminus spt(||V||) = \Omega_1 \cup \Omega_2 \cup \Omega_3$, where $\Omega_1$ is bounded by the loop $\gamma_2$, $\Omega_2$ is bounded by $\gamma_1$ and part of $\partial M$, and $\Omega_3$ is bounded by $\gamma_1\cup \gamma_2$ and a portion of the boundary of $M$. The sweepout will be constructed in four steps that we describe next.

Firstly, as in the proof of Proposition \ref{varrer-regioes}, there exists a sweepout $\{\sigma_t\}$, $t\in [0, \frac{1}{4}]$, of $\Omega_1$ such that $\sigma_0$ is a point in this region, $\sigma_{1/4} = \gamma_2$, and the lengths of all $\sigma_t$ are bounded by the length of $\gamma_2$.

For the second step, let $\gamma_2 : [0,\ell] \rightarrow M$ be a parametrization of the loop by arc length, with $\gamma_2(0) = \gamma_2(\ell) = p$. Suppose that the velocity vector of $\gamma_1$ emanating from $p$ is contained in the angle determined by $-\gamma_2^{\prime}(\ell)$ and a vector of $T_p(\partial M)$, such that this angle does not contain $\gamma_2^{\prime}(0)$. The stationarity of $V$ implies that the angle between $\gamma_2^{\prime}(0)$ and $T_p(\partial M)$ that does not contain the velocity vector of $\gamma_1$ is shorter than $\frac{\pi}{2}$. Then, consider $\{\sigma_t\}$, $t \in [\frac{1}{4}, \frac{1}{2}]$, which deforms $\gamma_2$ without increasing length to a curve $\sigma_{1/2}$ which is the sum of $\gamma_2$ restricted to $[\varepsilon, \ell]$ and the shortest geodesic from $\gamma_2(\varepsilon)$ to $\partial M$, for some small $\varepsilon$. Observe that this shortest segment from $\gamma_2(\varepsilon)$ to $\partial M$ is contained in $\Omega_3\cup\{\gamma_2(\varepsilon)\}$.   

In the third step, $\{\sigma_t\}$, $t \in [\frac{1}{2}, \frac{3}{4}]$, the curve $\sigma_{1/2}$ is deformed to the union $\sigma_{3/4} = \sigma_{1/2} + \gamma_1$. It suffices to add $\sigma_{1/2}$ to each curve of a sweepout of $\Omega_2$ by curves of lengths bounded by the length of $\gamma_1$.

Finally, observe that Lemma \ref{varrer-convexos} can be applied to the region of $M$ that remains to be swept. This region coincides with $\Omega_3$ with a small region removed. The removed region has a shape that resembles a triangle with vertices at $p$, $\gamma_2(\varepsilon)$ and at the end point of $\sigma_{1/2}$. After this region is removed from $\Omega_3$, the result is the only connected region whose boundary is $\sigma_{3/4}$. Therefore, one can obtain a family $\{\sigma_t\}$, $t \in [\frac{3}{4}, 1]$, deforming the curve $\sigma_{3/4}$ to a constant curve $\sigma_1$. Note that, it follows from the construction, and the propositions that were applied in the process, that these curves satisfy $L(\sigma_t) \leq L(\sigma_{3/4}) < ||V||(M)$, for all $t\in [0,1]$. 

In the case that the two segments of $V$ are type-(ii) loops at $p$, the construction is completely similar. It suffices to replace $\Omega_2$ by the region bounded by the second loop. 

\end{proof}

\subsection{Corollary}\label{interior}
{\em Let $V$ be as above with a singular point is of type ($J_i$); an interior junction such that $\theta^1(V,p) \geq 3$. Then, $M$ can be swept out by curves of lengths bounded by a number strictly lower than the mass of $V$.}

\begin{proof}
Consider the set $\Lambda \subset T_pM$ of unit vectors which are velocities of some geodesic segment emanating from $p$. In the case of a loop, there will be a pair of different vectors in $\Lambda$ associated with this curve.  If a given geodesic segment has multiplicity two, consider two copies of its velocity vectors in $\Lambda$, in such a way that the number of vectors of $\Lambda$ is exactly twice $\theta^1(V,p)$.

Recall that Proposition \ref{W-properties} implies that the support of $V$ intersects $\partial M$. In particular, it contains at least one type-(i) geodesic segment. Any pair of type-(i) segments determine exactly one angle in $T_pM$ whose measure is less than $\pi$. Otherwise, $V$ would contain a simple free boundary geodesic. In particular, if more than one type-(i) geodesic segment has multiplicity two, item (b) of Proposition \ref{Prop-multiplicities-angles} would give the desired result.

From now on, suppose that at most one type-(i) geodesic segment has multiplicity two. Fix a vector $v_1 \in \Lambda$, such that the geodesic emanating from $p$ with velocity $v_1$ is of type-(i). Choose an orientation for $T_pM$, and enumerate the other vectors of $\Lambda$ in the order of the positive angles that they make with $v_1$, in the chosen orientation. If some geodesic segment has multiplicity two, the repeated vectors are consecutive in this enumeration.

If the pair of different vectors in $\Lambda$ associated with a fixed loop are not consecutive in the list, then a geodesic segment starting at some other vector of $\Lambda$ is trapped in the loop. Therefore, there is a loop inside the first loop, which could coincide with the first one. In this case, item (a) of Proposition \ref{Prop-multiplicities-angles} implies the desired property.
 
In what follows, assume that the two different vectors in $\Lambda$ associated with any loop are consecutive in the enumeration: $v_i$ and $v_{i+1}$. Let us write $\theta(v_i, v_j)$ to represent the angle containing all the unit vectors in $T_pM$ whose positive angle with $v_i$ is between zero and the positive angle from $v_i$ to $v_j$. Next, it is shown that one of the following properties hold:

\begin{enumerate}
    \item[(I)] there exists a loop $\omega$ in $V$, with velocities at $p$ given by $v_i$ and $v_{i+1}$, whose angles $\theta(v_1, v_i)$ and $\theta(v_{i+1}, v_1)$ measure less than $\pi$;
    
    \item[(II)] there exist two different geodesic segments $\omega$ and $\eta$ of $V$, such that the first (possibly the unique) velocity vector of $\eta$ is the successor of the last (possibly the unique) velocity of $\omega$ in the enumeration, such that the angles $\theta(v_1, v_i)$ and $\theta(v_j, v_1)$ measure less than $\pi$, where $v_i$ is the first velocity vector of $\omega$ and $v_j$ is the last of $\eta$.     
\end{enumerate}
 
Indeed, choose $\omega$ to be the geodesic segment such that its first velocity vector at $p$, $v_i$, satisfies $\theta(v_1, v_i) < \pi$, and $i$ is the largest index with this property in the enumeration. If $\omega$ is of type-(i), then the starionarity implies that $v_i$ is not the last vector in $\Lambda$. In this case, let $\eta$ be the segment whose first velocity at $p$ is $v_{i+1}$. It follows that $\theta(v_1, v_{i+1})\geq \pi$. If equality holds, then $\eta$ is a loop, as $V$ does not contain simple free boundary geodesics. Then, property (II) holds. Suppose now that $\omega$ is a loop. In this case, if $\theta(v_1, v_{i+1}) \leq \pi$, then the same argument above shows that one can guarantee that (II) holds true. Then only remaining case is exactly the validity of (I) for the choice of $\omega$ that was made.

Let us consider case (I). Observe that the assumption that the multiplicity of the loop is one, and the fact that the density $\theta^1(V,p)$ of $V$ at $p$ is at least three, imply that there are at least three vectors in $\Lambda$ other than $v_1$, and the velocities $v_i$ and $v_{i+1}$ of $\omega$. Then, at least two of them belong either to $\theta(v_1, v_i)$ or $\theta(v_{i+1}, v_1)$. These two are either a pair of type-(i) segments, or a loop inside and angle measuring less than $\pi$. In any situation, part (b) of Proposition \ref{Prop-multiplicities-angles} ends the argument.

Case (II) is divided into three particular cases, depending on the types of $\omega$ and $\eta$. Assume first that $\omega$ and $\eta$ are both of type-(i). As before, the density assumption implies that at least three other vectors of $\Lambda$ belong to one of the angles $\theta(v_1, v_i)$ or $\theta(v_j, v_1)$, both measuring less than $\pi$. Part (b) of Proposition \ref{Prop-multiplicities-angles} can be applied as before to conclude the argument. Assume now that $\omega$ is a loop and $\eta$ is of type-(i). Then, at least two other vectors of $\Lambda$ belong to one of the angles $\theta(v_1, v_i)$ or $\theta(v_j, v_1)$. If both of these extra vectors belong to the same angle, part (b) of Proposition \ref{Prop-multiplicities-angles} can be applied. Assume now that $2<i<j<k$, for some $k \leq 2\cdot \theta^1(V, p)$. Then, the shortest angle between $v_2$ and $v_j$, which measures less than $\pi$, contains either the velocities of the loop $\omega$, or the pair $v_1$ and $v_k$. In any case, one can invoke part (b) of Proposition \ref{Prop-multiplicities-angles}. The case in which $\eta$ is a loop and $\omega$ is of type-(i) is completely similar.

Finally, consider the case in which (II) holds for two loops $\omega$ and $\eta$. Since the loops have multiplicity one in $V$, the regions $\Omega_{\omega}$ and $\Omega_{\eta}$ bounded by these curves belong to the same side of the decomposition $\mathcal{I}\cup\mathcal{J}$ obtained in Lemma \ref{W-properties} with $W = V$. Let $\Omega_1$ be the component of $M\setminus spt(||V||)$ which bounds from the exterior the loops $\omega$ and $\eta$. The existence of such a region follows from the fact that these loops are consecutive in the enumeration of segments emanating from $p$ that was considered before. Moreover, if $\Omega_{\omega}$ and $\Omega_{\eta}$ belong to $\mathcal{I}$, then $\Omega_1$ belongs to $\mathcal{J}$. The desired sweepout is constructed applying methods that are similar to those used at the end of the proof of Corollary \ref{bdry}. First one sweeps $\Omega_{\omega}\cup \Omega_{\eta}$ out as in the proof of Proposition \ref{varrer-regioes}. The final curve being $\omega+ \eta$. Next, deform this curve decreasing length, by a one-parameter family, to a curve which is the boundary of $\Omega_{\omega}\cup\Omega_{\eta}\cup R$, where $R$ is a small triangle determined by geodesics connecting $p$ and one point of each loop $\omega$ and $\eta$. This can be thought as deforming the union of consecutive loops inside the angle bounded by the second velocity vector of $\omega$ and the first of $\eta$. One replaces the union of small geodesic segments of $\omega$ and $\eta$ with one extremity at $p$ by the minimizing geodesic connecting the second extremities of each segment. Next, sweep all the other components of $\mathcal{I}$ out. As in the proof of Proposition \ref{varrer-regioes}, the sweepout constructed so far is combined at the final curve to another sweepout obtained from the regions in $\mathcal{J}$. The difference is that region $\Omega_1$, exterior to the loops $\omega$ and $\eta$, is replaced by the slightly smaller region $\Omega_1 \setminus R$. This is possible because the latter still satisfies the convexity properties necessary for the application of Lemma \ref{varrer-convexos}. The interpolation can be performed in such a way that lengths increase less than the gap created by the substitution of the two small geodesic segments with vertex at $p$ by a single minimizing segment.   
\end{proof}

\subsubsection*{Remark:} The sweepouts considered at the end of the proofs of Corollaries \ref{bdry} e \ref{interior} generalize those considered by Pitts in arguments similar to those of Lemma 2.3 of \cite{CalCao}. The generalization is in the sense that in the present work one does not know if both exterior angles determined by the pair of loops are convex, and multiplicity different from one must be considered.

In what follows, we consider the case of a geodesic network in $M$ which has exactly two geodesic segments meeting at a point $p \in \partial M$. This is precisely the case of a junction of type ($J_{\ell}$), as introduced in \ref{admissible}. 

Let $(M^2,g)$ be isometrically embedded in a closed Riemannian surface $N$. Let $v_1, v_2  \in T_p M$ be unit vectors pointing inside $M$, $v_1\neq v_2$, and $\gamma_i$ be the geodesics of $N$ with $\gamma_i(0)=p$ and $\gamma_i^{\prime}(0)=v_i$, $i=1, 2$. Consider the varifold $V$ defined as the sum of the varifolds induced by the restrictions $\gamma_i|_{[0,T(i)]}$ of $\gamma_i$ to the largest positive interval containing zero and with $\gamma_i([0,T(i)])\subset M$. Assume that $V$ is stationary. This implies that $\gamma_i^{\prime}(T(i))$ is perpendicular to $\partial M$ at $\gamma_i(T(i))$, for $i=1,2$, and the vector  $\gamma_1^{\prime}(0) + \gamma_2^{\prime}(0) = v_1+v_2$ is perpendicular to $\partial M$ at $p$.

\subsection{Proposition}\label{V-figure}
{\em Let $V$ be a stationary geodesic network composed of two segments meeting at a boundary singular point $p\in \partial M$, as considered above. Then, $M$ can be swept out by curves of lengths bounded by a number strictly lower than the mass of $V$.}

The argument applied in the proof is motivated by the constructions performed in the proof of Claim 1 related to Theorem D of \cite{CalCao}. It is also similar to the argument sketched in the proof of Proposition \ref{varrer-convexos}. 

\begin{proof}
Consider the notation introduced in the paragraph immediately above the statement of the present proposition. Let $n_i$ be a choice of unit normal vector field along $\gamma_i$ such that $n_i(t) \in T_{\gamma_i(t)}M$ for all $t$, and that the vectors $n_1(0)$ and $-n_2(0)$ point inside $M$. The stationarity condition implies that the sum of the normal vectors satisfy $n_1(0)+n_2(0)\in T_p(\partial M)$.

The signed distance function to $\partial M$ in $N$ is denoted by $d$. This function is positive inside $M$, and is the negative of the distance to $\partial M$ in $N\setminus M$. This function is smooth in a neighborhood of $\partial M$, and its gradient $\nabla d$ is orthogonal to the equidistant surfaces (levels of $d$). Define, for $i= 1,2$, 
$$X_i(t) = n_i(t) - \frac{\langle n_i(t), \nabla d\rangle}{\langle \gamma_i^{\prime}(t), \nabla d\rangle } \gamma_i^{\prime}(t),$$
where $\nabla d$ is computed at $\gamma_i(t)$. Observe that $X_i(t)$ is perpendicular to $\nabla d(\gamma_i(t))$, and $X_1(0) = X_2(0)$. Indeed, $X_i(0)$ is the projection of the vector $n_i(0)$ onto $T_p(\partial M)$ parallel to the line spanned by $v_i$. Moreover, since $v_1 \neq v_2$ and $(v_1+v_2) \perp T_p(\partial M)$, it follows that $X_i(t)\neq 0$, for small $t$.

Consider the vector fields $P_{12}$ and $P_{21}$ defined as the parallel transports of $v_1$ along $\gamma_2$ and $v_2$ along $\gamma_1$, respectively. Here, these vector fields are considered even at points of $\gamma_i$ outside $M$. For $q$ in a neighborhood of $p$ inside $N$, there exists a unique  $\pi_i(q) \in \gamma_i$, such that $exp_{\pi_1(q)}(s_1 P_{21}) = q$ and $exp_{\pi_2(q)}(s_2 P_{12}) = q$, for some small numbers $s_1$ and $s_2$, where both vector fields are applied at the corresponding points $\pi_i(q)$. Moreover, the maps $\pi_i(q)$ are smooth, and $\pi_1(\gamma_2(t)) = p = \pi_2(\gamma_1(t))$. Observe that it might happen that $\pi_i(q)$ does not belong to $M$. Consider the function
$$
\varphi(q) = |X_1(\pi_1(q))| + |X_2(\pi_2(q))| - |X_1(0)|,
$$
where $X_i(\pi_i(q)) = X_i(t)$ with $\gamma_i(t) = \pi_i(q)$. Since $X_i(t)\neq 0$ for small $t$, it follows that this function is smooth in a neighborhood of $p$. Since $X_1(0) = X_2(0)$, one has $\varphi(\gamma_i(t)) = |X_i(t)|$, for small $t$.

Let $Z = Z(q)$ be the unit vector field defined in a neighborhood of $\partial M$ such that $\langle Z, \nabla d\rangle = 0$ and $X_i(0) = c Z(p)$, for some positive constant $c$.

Finally, set $X(q) = \varphi(q)Z(q)$, for $q$ in a small neighborhood of $p$. This vector field extends $X_i(t)$ along $\gamma_i$, $i = 1,2$, to a smooth vector field which is tangential to $\partial M$. Using smooth extensions of the vector fields $n_i$ to neighborhoods of $\gamma_i$, and a partition of unity argument one can construct a vector field $Y$ supported in a neighborhood of the support of $V$ satisfying:
\begin{enumerate}
    \item[(i)] $Y(q) \in T_q(\partial M)$, for $q \in \partial M$.
    
    \item[(ii)] The component of $Y(\gamma_i(t))$ orthogonal to $\gamma_i$ has unit length.
\end{enumerate}
The vector field $Y$ coincides with $X$ in a neighborhood of $p$, and with $n_i$ along the restriction of $\gamma_i$ to some interval of the form $[t_i, T(i)]$, with $0< t_i < T(i)$. In particular, $Y(\gamma_i(T(i))) = n_i(T(i))$.

Let $\{F(\cdot, t)\}$ denote the flow of the vector field $Y$, and consider the push-forward $V_t = F(\cdot,t)_{\#}V$ of $V$ by this one-parameter family of maps. Property (i) implies that $spt(V_t)\subset M$. The second variation formula and (ii) imply
$$
\frac{d^2}{dt^2} L(V_t) \big|_{t=0} = - \int_{0}^{T(1)} K(\gamma_1(t))dt - \int_{0}^{T(2)} K(\gamma_2(t))dt + \sum \langle \nabla_Y Y, \nu_i(t)\rangle,
$$
where the sum in the formula is over $i = 1$ and $2$, and $t = 0$ and $t = T(i)$, and $\nu_i(t)$ represents the exterior unit co-normal of $\gamma_i$ at the extremity $\gamma_i(t)$.

At the regular boundary points $\gamma_i(T(i))$, the vector $\nu_i(T(i))$ is orthogonal to $\partial M$. This implies that $\langle \nabla_Y Y, \nu_i(T(i))\rangle = - A_{\partial M}(Y,Y)$. The stationarity condition at $p$ implies that $\nu_1(0)+\nu_2(0) = - v_1 - v_2$ is also perpendicular to $\partial M$ and outward pointing. Therefore,
$$
\langle \nabla_Y Y, \nu_1(0) + \nu_2(0)\rangle = - A_{\partial M}(Y(p),Y(p))|\nu_1(0)+\nu_2(0)|.
$$
Putting all this information together, and combining with the convexity and positive curvature assumption, it follows that $(d^2/dt^2)|_{t=0} L(V_t) <0$.

The remaining of the construction is similar to that applied in Steps 2 and 3 in the proof of Claim 1 related to Theorem D of \cite{CalCao}. We include the steps of the adapted construction for completeness.  
Fix $\varepsilon>0$ small, and parametrize the supports of $V$ and $V_{\varepsilon}$ respectively by maps $\gamma$ and $\gamma_{\varepsilon}$ defined on $[-1,1]$ such that $\gamma(0)=p$, $\gamma(-1) = \gamma_1(T(1))$, $\gamma(1) = \gamma_2(T(2))$,  
$$
\gamma([-1,0]) = \gamma_1([0, T(1)]) \text{ and } \gamma([0,1]) = \gamma_2([0, T(2)]),
$$
and, similarly for the deformed curve, assume that $\gamma_{\varepsilon}(0)=F(p, \varepsilon)$, $\gamma_{\varepsilon}(-1) = F(\gamma_1(T(1)), \varepsilon)$, $\gamma_{\varepsilon}(1) = F(\gamma_2(T(2)), \varepsilon)$, with 
$$
\gamma_{\varepsilon}([-1,0]) = F(\gamma_1([0, T(1)]), \varepsilon) \text{ and } \gamma_{\varepsilon}([0,1]) = F(\gamma_2([0, T(2)]), \varepsilon).
$$
For simplicity, let $\eta^+ = \gamma_{\varepsilon}|_{[\varepsilon,1]}$ and $\eta^- = \gamma_{\varepsilon}|_{[-1,-\varepsilon]}$. Next, consider
$$
\eta_{\varepsilon} = \eta^+ + \tau(\gamma_{\varepsilon}(-\varepsilon), \gamma_{\varepsilon}(\varepsilon)) + \eta^-,
$$
and
$$
\eta_{-\varepsilon} = \eta^+ + \tau(\gamma_{\varepsilon}(\varepsilon), \partial M) + \eta^- + \tau(\gamma_{\varepsilon}(-\varepsilon), \partial M),
$$
where $\tau(x,y)$ represents the shortest geodesic of $M$ from point $x$ to point $y$, and $\tau(x,\partial M)$ represents the shortest geodesic of $M$ from point $x$ to $\partial M$. 

There exists a homotopy $\{\eta_{t}\}$, $-\varepsilon\leq t \leq \varepsilon$, from $\eta_{-\varepsilon}$ to $\eta_{\varepsilon}$ by curves which satisfy $L(\eta_t) \leq L(V_{\varepsilon}) < ||V||(M)$.

One can perform a similar cut-and-paste construction on $\gamma$. Letting $\gamma^- = \gamma|_{[-1, -\varepsilon]}$ and $\gamma^+ = \gamma|_{[\varepsilon, 1]}$  denote restriction of $\gamma$, define
$$
\eta_{2\varepsilon} = \gamma^+ + \tau(\gamma(-\varepsilon), \gamma(\varepsilon)) + \gamma^-
$$
and
$$
\eta_{-2\varepsilon} = \gamma^+ + \tau(\gamma(\varepsilon), \partial M) + \gamma^- + \tau(\gamma(-\varepsilon), \partial M).
$$
There exists a homotopy $\{\eta_t\}$, $t\in [-2\varepsilon, -\varepsilon] \cup [\varepsilon, 2\varepsilon]$ connecting $\eta_{-2\varepsilon}$ and $\eta_{-\varepsilon}$, and connecting $\eta_{\varepsilon}$ and $\eta_{2\varepsilon}$. Moreover, this can be obtained using curves whose length satisfy $L(\eta_t) \leq ||V||(M) - \delta$, for some $\delta>0$. 

The end of the construction consists of sweeping out the convex regions of $M$ bounded by the broken geodesics $\eta_{-2\varepsilon}$ and $\eta_{2\varepsilon}$ using Lemma \ref{varrer-convexos}.
\end{proof}

\section{Proof of Theorem \ref{thm-main} and two inequalities}\label{sec-proof}

In this brief section, the results used in the proof of the main theorem of the paper are combined. The proofs of two inequalities stated in the introduction of the paper are also included.

\begin{proof}[Proof of Theorem \ref{thm-main}]
It follows from Theorem \ref{thm-single-jct} that there exists a free boundary stationary geodesic network $V$ with at most one junction, denoted by $p$ whenever it exists, such that $\omega(M) = ||V||(M)$.

Consider the case where $p \in int(M)$ is a singular point at which the density satisfies $\theta^1(V,p)=2$. Theorem \ref{thm-single-jct} implies that $p$ is a crossing point of two geodesic segments with multiplicity one. In particular, $V$ would be the crossing point of two local geodesics. These could be part of a single global geodesic with a self-intersection. In any case, $V$ would be either the union of two simple free boundary geodesics (possibly two copies of the same curve), or a free boundary geodesic with a self-intersection in the interior. Proposition \ref{varrer-regioes} would imply that $V$ does not realize the width.

The case of an interior junction $p \in int(M)$ at which $\theta^1(V,p)\geq 3$ is ruled out by Corollary \ref{interior}. The case of $p\in \partial M$ at which $\theta^1(V,p)\geq 1.5$ is ruled out by Corollary \ref{bdry}. The case $p \in \partial M$, $\theta^1(V,p)= 1$, and $V$ is not a geodesic loop is ruled out by Proposition \ref{V-figure}. This finishes the argument.
\end{proof}

\begin{proof}[Proof of $\omega(M) < L(\partial M)$]
Let $\beta : [0, T]\rightarrow \partial M$ be a parametrization by arc length of the boundary curve, where $T = L(\partial M)$. For large $n$, the sequence $\{p_i\}$, $1\leq i \leq n$, of points $p_i = \beta(\frac{iT}{n})$ is such that there exists a unique minimizing geodesic from $p_i$ to $p_{i+1}$, for $1\leq i \leq n-1$, and from $p_n$ to $p_1$. Let $\rho$ be this broken geodesic. Observe that $M\setminus \rho$ is made of $n+1$ convex regions. The bigger region, which is denoted by $U_0$, has boundary $\rho$. The other components, $\{U_i\}$, $1\leq i \leq n$, are bounded by one of the minimizing geodesic segments and a short arc of $\partial M$. The regions $U_i$, $i\neq 0$, can be swept out using Lemma \ref{varrer-convexos}. The region $U_0$ can be swept out using the remark after Lemma \ref{varrer-convexos}. The sweepout of $M$ obtained by the concatenation of the sweepouts of each $U_i$ is such that the supremum of the lengths of all curves is at most $L(\rho)$. In addition, $L(\rho) < L(\partial M)$.
\end{proof}

\begin{proof}[Proof of $\omega_4(\mathbb{D}) < \omega_5(\mathbb{D})$]
Here $\mathbb{D}$ denotes the flat disk of unit radius in $\mathbb{R}^2$. Consider the decomposition $\mathbb{D} = \cup_{i=1}^5 S_i$ of $\mathbb{D}$ by five congruent circular sectors of central angle $\frac{2\pi}{5}$. Approximate each $S_i$ by compact strictly convex regions $K(i,j)$, $j\in \mathbb{N}$, contained in the relative interior of $S_i$. Theorem \ref{thm-main} implies that $\omega(K(i,j))$ is realized by the length of a free boundary segment $V(i,j)$. As in the case of triangles, these regions can be swept out by parallel segments with lengths bounded by $h = \sin(\frac{2\pi}{5})$, see Figure \ref{Fig-2}. Therefore, the segments $V(i,j)$ do not converge to a radius of $S_i$, or to a chord connecting the endpoints of the arc of circle in Figure \ref{Fig-2}, since these segments have length bigger than $h$. In conclusion, the $V(i,j)$ converge, as $j\rightarrow \infty$, precisely to the highlighted segment of Figure \ref{Fig-2} whose length is equal to $h$. The Lusternik-Schnirelmann inequality, see Liokumovich, Marques, and Neves \cite{LioMarNev}, Section 8.3 in Gromov \cite{Gro}, and Section 3 in Guth \cite{Gut}, tells us that
$$
\omega_5(\mathbb{D}) \geq \sum_{i=1}^5 \omega_1(K(i,j)). 
$$
Letting $j$ tend to infinity, one obtains that $\omega_5(\mathbb{D}) \geq 5h > 4 = \omega_4(\mathbb{D})$.

\begin{figure}[htp]
    \centering
    \includegraphics[width=4cm]{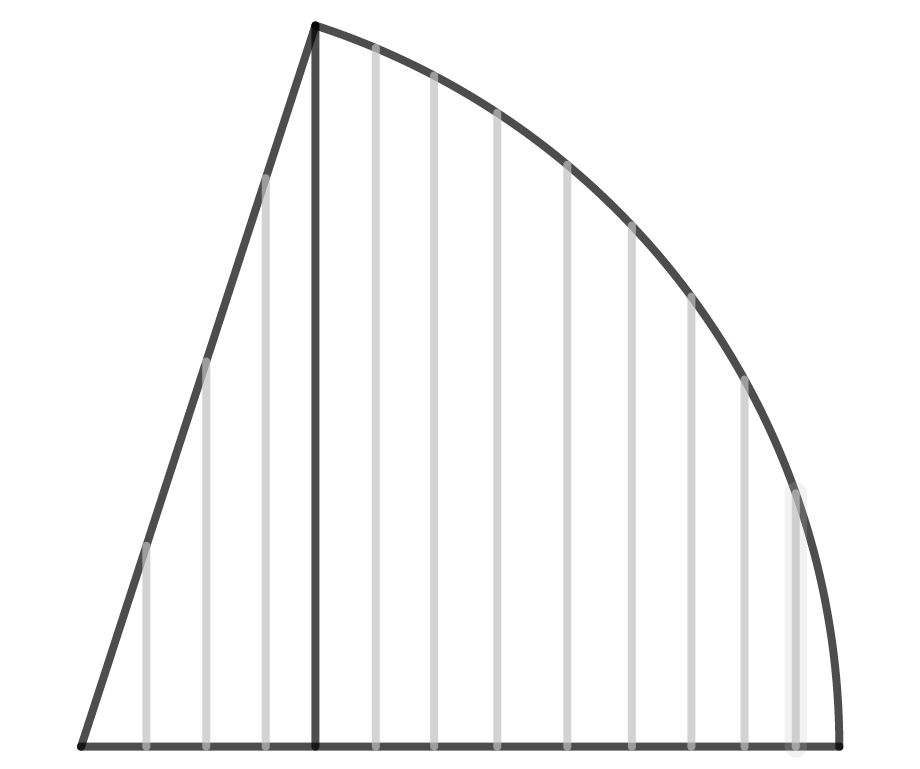}
    \caption{A sweepout of the circular sector by parallel segments.}
    \label{Fig-2}
\end{figure}
\end{proof}

\bibliographystyle{amsbook}

\end{document}